\documentclass[a4paper,11pt]{article}

\usepackage[T1]{fontenc}
\usepackage{lmodern}

\usepackage{dsfont}

\usepackage[latin1]{inputenc}
\usepackage{amsmath}
\usepackage{amsthm}
\usepackage{amssymb}
\usepackage{mathrsfs}
\usepackage{graphicx}
\usepackage[all]{xy}
\usepackage{hyperref}

\usepackage{makeidx}

\usepackage{stmaryrd}
\usepackage{caption}

\usepackage{abstract} 

\newtheorem{thm}{Theorem}

\theoremstyle{definition}

\usepackage{chngcntr}
\counterwithout{paragraph}{subsubsection}

\def\rquotient#1#2{%
	\makeatletter
	\raise.3ex\hbox{$#1$}/\lower.3ex\hbox{$#2$}%
	\makeatother
}	

\makeatletter
\newcommand{\subjclass}[2][2010]{%
	\let\@oldtitle\@title%
	\gdef\@title{\@oldtitle\footnotetext{#1 \emph{Mathematics subject classification.} #2}}%
}
\newcommand{\keywords}[1]{%
	\let\@@oldtitle\@title%
	\gdef\@title{\@@oldtitle\footnotetext{\emph{Key words and phrases.} #1.}}%
}
\makeatother

\newcommand{\Address}{{% additional braces for segregating \footnotesize
		\bigskip
		\small
		
		\textsc{University of Montpellier\\ 
Institut Math\'ematiques Alexander Grothendieck\\
Place Eug\`ene Bataillon\\
34090 Montpellier (France)}\par\nopagebreak
		\textit{E-mail address}: \texttt{anthony.genevois@umontpellier.fr}
		
}}

\makeindex

\title{Why CAT(0) cube complexes should be replaced with median graphs}
\date{\today}
\author{Anthony Genevois}

\begin{document}

\maketitle

\begin{abstract}
In this note, we discuss and motivate the use of the terminology ``median graphs'' in place of ``CAT(0) cube complexes'' in geometric group theory.
\end{abstract}

%\tableofcontents

\vspace{1cm}
A \emph{CAT(0) cube complex} refers to a cube complex (i.e.\ a cellular complex obtained by gluing cubes of various dimensions along faces) whose length metric obtained by extending the Euclidean metrics on its cubes (by convention with sides of length one) is CAT(0). They first appear in \cite{MR0919829}, where they are introduced as a convenient source of CAT(0) spaces. Indeed, while it is in general a difficult task to determine whether or not a given geodesic metric space is CAT(0), CAT(0) cube complexes can be identified thanks to a simple local criterion:

\begin{thm}[\cite{MR0919829, MR3029427}]\label{thm:Local}
A cube complex is CAT(0) if and only if it is simply connected and all its vertices have flag simplicial links.
\end{thm}

Roughly speaking, the \emph{link} of a vertex corresponds to a small sphere around it endowed with the complex structure induced by the cubical structure of the whole complex. Thus, vertices (resp. edges, triangles, $n$-simplices) in the link correspond to edges (resp. squares, $3$-cubes, $(n+1)$-cubes) in the cube complex. A simplicial complex is \emph{flag} if any pairwise adjacent vertices span a simplex. A typical example where the local criterion provided by Theorem~\ref{thm:Local} fails is an $n$-cube minus its $n$-cell (so, topologically, an $(n-1)$-sphere). 

An important feature of CAT(0) cube complexes is highlighted in \cite{MR1347406}: CAT(0) cube complexes admit natural ``codimension-one separating subspaces'', called \emph{hyperplanes}. The graph-metric on the one-skeleton of a CAT(0) cube complex, which is quasi-isometric to the CAT(0) metric in the finite-dimensional case, turns out to be tightly connected to hyperplanes. For instance, a combinatorial path between vertices in the one-skeleton is a geodesic (with respect to the graph-metric) if and only if it crosses each hyperplane at most once. As a consequence, the graph-distance between two vertices coincides with the number of hyperplanes separating them. Since then, the combinatorics of hyperplanes, and consequently the graph-metric, has been a central tool in the study of CAT(0) cube complexes.

Interestingly, the one-skeletons of CAT(0) cube complexes coincide with a well-known family of graphs: \emph{median graphs}. Formally introduced in the 1960s-70s \cite{MR0125807, MR0286705, MR522746}, but with even older roots motivated by lattice theory, median graphs have been extensively investigated in metric graph theory. (See for instance \cite{MR2405677, MR2798499} and references therein.)

\begin{thm}[\cite{MR1663779, MR1748966, Roller}]
The one-skeleton of a CAT(0) cube complex is a median graph. Conversely, the cube-completion of a median graph (i.e.\ the cube complex obtained by filling with cubes all the subgraphs isomorphic to one-skeletons of cubes) is a CAT(0) cube complex.
\end{thm}

Recall that a graph is \emph{median} if, for any three vertices $x_1,x_2,x_3$, there exists a unique vertex $m$ (the \emph{median point}) that lies at the intersection of three geodesics connecting $x_1,x_2,x_3$, i.e.\ satisfying $d(x_i,x_j) = d(x_i,m)+d(m,x_j)$ for all distinct $i,j \in \{1,2,3\}$. 

Despite the fact that CAT(0) cube complexes and median graphs essentially define the same objects, and the fact that CAT(0) cube complexes are now mainly studied through their one-skeletons, the terminology used in geometric group theory remains ``CAT(0) cube complexes'' in the literature. In this note, we would like to motivate the use of ``median graphs'' instead.

\paragraph{CAT(0) geometry is not necessary.} Very few articles actually use CAT(0) geometry in an essential way to prove results about CAT(0) cube complexes. Instead, they mainly use the graph-metric and the combinatorics of hyperplanes. Nevertheless, before replacing ``CAT(0) cube complexes'' with ``median graphs'', at the very least we have to verify that the central results about the structure of CAT(0) cube complexes and about group actions used in geometric group theory can be easily and naturally proved by using the point of view of median graphs. This is done, for instance, in \cite{Book}. 

 In order to justify the replacement of CAT(0) cube complexes with median graphs, we do not have to prove that CAT(0) geometry will never be relevant. Actually, there may (and probably do) exist a few specific cases where CAT(0) geometry on cube complexes is superior from some perspective\footnote{Personally, I am not aware of a single theorem about CAT(0) cube complexes that cannot be naturally proved thanks to median graphs.}. But it suffices to show, as tried in this note and in \cite{Book}, that median geometry is more natural, efficient, and satisfying in most cases.

\paragraph{CAT(0) geometry not easily accessible.} Many group actions on CAT(0) cube complexes are obtained by cubulating a space with walls (e.g.\ Coxeter groups \cite{MR1983376}, small cancellation groups \cite{MR2053602}, hyperbolic $3$-manifold groups \cite{MR2931226}, one-relator groups with torsion \cite{MR3118410}); by taking the universal cover of a nonpositively cube complex (e.g.\ right-angled Artin groups \cite{MR1368655}, groups of alternating links \cite{MR1707529}, graph braid groups \cite{MR2701024}); by constructing directly the cube complex (e.g.\ graph products \cite{MR1389635}, diagram groups \cite{MR1978047}, Artin groups of infinite type \cite{MR3993762}, Cremona groups \cite{MR4340723}, Neretin groups \cite{Neretin}); or by finding a commensurating action (e.g.\ groups of local similarities \cite{MR2486801}, some topological full groups \cite{MR3377390}, Neretin groups \cite{Boudec}, birational groups \cite{MR4014636}, piecewise smooth homeomorphism groups \cite{MR4100128}, Grigorchuk groups \cite{GrigorCC}). In all these cases, computing the CAT(0) distance between two arbitrary points, or even vertices, is often out of reach. The graph-distance between two vertices, however, may be computable and even meaningful. 

For right-angled Artin groups, the median graphs on which they act are the Cayley graphs given by their canonical generating sets. Therefore, the distance between vertices coincides with the corresponding word-distance. For Coxeter groups, the median graphs are obtained by cubulating the canonical wallspace structures on their Cayley graphs. Since the distance in the median graph of two points coming from the wallspace coincides with the number of walls separating them, it follows that the distance in the median graph of two points coming from the Coxeter group agrees with the word-distance. 

As a more interesting example, consider the group $\mathrm{PDiff}(\mathbb{S}^1)$ of piecewise differentiable homeomorphisms $\mathbb{S}^1 \to \mathbb{S}^1$. We fix an orientation on $\mathbb{S}^1$, and, for all $g \in \mathrm{PDiff}(\mathbb{S}^1)$ and $x \in \mathbb{S}^1$, we denote by $g'(x^-)$ (resp. $g'(x^+)$) the left-derivative (resp. right-derivative) of $g$ at $x$. One can make $\mathrm{PDiff}(\mathbb{S}^1)$ act on $\mathbb{S}^1 \times (0,+ \infty)$ via
$$g \cdot (x,r) := \left( g(x), g'(x^-)rg'(x^+)^{-1} \right), \ g \in \mathrm{PDiff}(\mathbb{S}^1), (x,r) \in \mathbb{S}^1 \times (0,+\infty).$$
A key observation is that $\mathrm{PDiff}(\mathbb{S}^1)$ commensurates $S:= \mathbb{S}^1 \times \{1\}$, i.e.\ the symmetric difference $S \triangle gS$ is finite for every $g \in \mathrm{PDiff}(\mathbb{S}^1)$. This follows from the observation that
$$S \backslash g^{-1}S = \{ (x,1) \mid x \in \mathbb{S}^1 \text{ such that } g'(x^-) \neq g'(x^+) \}$$
for every $g \in \mathrm{PDiff}(\mathbb{S}^1)$. As a consequence, $\mathrm{PDiff}(\mathbb{S}^1)$ admits an action on a median graph $X$ such that, for some vertex $o \in X$, the equalities
$$d_X(o,g \cdot o ) = S \triangle gS = 2 \# \underset{=: \mathrm{Sing}(g)}{\underbrace{\{ \text{points at which $g$ is not differentiable} \} }}$$
hold for every $g \in \mathrm{PDiff}(\mathbb{S}^1)$. Because every isometry of a median graph with unbounded orbits admits (up to taking a convenient subdivision) a bi-infinite geodesic on which it acts as a translation, it follows that, for every $g \in \mathrm{PDiff}(\mathbb{S}^1)$, either $(\# \mathrm{Sing}(g^n))_{n \geq 1}$ is bounded or there exists a positive integer $K \geq 1$ such that
$$\# \mathrm{Sing}(g^n)= \frac{K}{2} \cdot n + O(1).$$
The constant $K$ corresponds to the translation length of $g$ in the median graph.

In the same spirit, consider, given a field $k$, the \emph{$d$-dimensional Cremona group $\mathrm{Cr}_k^d$ over $k$}, i.e.\ the group of birational transformations $\mathbb{P}_k^d \dashrightarrow \mathbb{P}_k^d$. In \cite{MR4340723}, a median graph on which $\mathrm{Cr}_k^d$ acts is explicitly constructed. The distance between vertices can be computed \cite[Lemma~4.11]{MR4340723}, and, combined with some median geometry, a linear lower bound on the asymptotic growth of $(\mathrm{deg}(f^n))_n$ can be deduced for non-pseudo-regularisable birational transformations $f : \mathbb{P}_k^d \dashrightarrow \mathbb{P}_k^d$. (Recall that $f$ has \emph{degree $k$} if one can write $f : [x_0,\ldots, x_d] \dashrightarrow [g_0, \ldots, g_d]$ for homogeneous polynomials $g_0, \ldots, g_d$ of degree $k$ without a non-constant common factor.) 

In the same way that distances with respect to the CAT(0) metric are usually not computable, the solutions known for various algorithmic problems in CAT(0) groups, such as the word and conjugacy problems, do not provide explicit procedures in practice. However, in some cubulable groups for which the structure of the median graph is tightly connected to the algebraic structure of groups (e.g.\ when the median graph is a Cayley graph), then it may be possible to exploit the median geometry in order to solve explicitly some algorithmic problems. This is illustrated, for instance, by the solution of the conjugacy problem in cactus groups in \cite{Cactus}.

\paragraph{Median geometry may be enlightening.} Using the point of view of median graphs may also lead to the introduction of new concepts, which would not have been accessible from CAT(0) geometry. This can be illustrated by the \emph{Roller boundary}, which can be described as follows \cite{MR4071367}. Let $X$ be a countable median graph. Fix a basepoint $o \in X$. The Roller boundary $\mathfrak{R}X$ is the graph
\begin{itemize}
	\item whose vertices are classes of geodesic rays starting from $o$, where two rays are equivalent whenever they cross exactly the same hyperplanes;
	\item whose edges connect two (classes of) rays whenever there exists exactly one hyperplane that is crossed by one but not by the other.
\end{itemize}
The Roller boundary is usually not connected, but each connected component turns out to be a median graph. Thus, one gets the same structure on the space and on its boundary\footnote{Even better, the components are organised along a median graph! More precisely, let $X$ be a median graph. Consider the Hasse diagram $\mathfrak{C}X$ of the poset whose elements are the components of the Roller completion $\overline{X}:= X \cup \mathfrak{R}X$, and whose order is defined as follows. Given two components $Y,Z$ of $\overline{X}$, $Y \leq Z$ if $Y$ is contained in the image of the Roller completion $\overline{Z}$ of $Z$ in $\overline{X}$. Then every connected component of $\mathfrak{C}X$ is a median graph. If $X$ has finite cubical dimension, $\mathfrak{C}X$ is connected and bounded.}, contrasting with visual boundaries of CAT(0) spaces. This property may be useful in some inductive arguments: if a group acts on a median graph $X$ (of finite cubical dimension) and stabilises a component $Y$ of $\mathfrak{R}X$, then it acts on a new median graph $Y$ (of smaller cubical dimension). See for instance \cite{MR4062290} for an application of this idea. 

Let us illustrate how the Roller boundary may be useful, and conceptually satisfying, by sketching a proof of the Tits alternative extracting from \cite{MR2827012}\footnote{Here, we follow \cite{Book}. In particular, we avoid the use of the visual and simplicial boundaries. See also \cite{MR3959849} for a similar approach applied to median spaces of finite rank.}. So let $G$ be a group acting on a median graph $X$ with finite cubical dimension. We claim that $G$ either contains a non-abelian free subgroup or is virtually (locally $X$-elliptic)-by-(free abelian).
\begin{itemize}
	\item A hyperplane $J$ is \emph{$G$-flippable} if there exists $g \in G$ such that $gD \subset D^c$ for some halfspace $D$ delimited by $J$. By an easy ping-pong argument, we know that, if $G$ acts on a median graph with a facing triple (i.e.\ three hyperplanes such that no one separates the other two) of $G$-flippable hyperplanes, then $G$ contains a non-abelian free subgroup.
	\item If some hyperplane $J$ is not $G$-flippable, then, fixing a halfspace $X$ delimited by $J$, the halfspaces $gD$ with $g \in G$ pairwise intersect. The intersection $\bigcap_{g \in G} gD$ either is non-empty, providing a $G$-invariant convex subgraph; or accumulates to the Roller boundary, providing a $G$-invariant component of $\mathfrak{R}X$. By iterating the argument, one finds a $G$-invariant convex subgraph $Y$ in $\overline{X}:= X \cup \mathfrak{R}X$ all of whose hyperplanes are $G$-flippable.
	\item If $Y$ contains a facing triple, then we know from the first item that $G$ contains a non-abelian free subgroup. Otherwise, we can show that $Y$ must have only finitely components in its Roller boundary (like the infinite grid $\mathbb{Z}^n$), which implies that $G$ virtually fixes a point $\xi$ in $\mathfrak{R}X$. 
	\item Let $Z_1, \ldots, Z_k$ denote the components of $\mathfrak{R}X$ containing $\xi$ or containing $\xi$ in their boundaries. We can prove that there are only finitely many such components. Because $G$ virtually fixes $\xi$, it contains a finite-index subgroup $H$ stabilising each $Z_i$. For every $1 \leq i \leq k$, we can define a morphism $H \to \mathbb{Z}$ by
$$\mathfrak{h}_i : h \mapsto |\mathcal{W}_i\backslash g \mathcal{W}_i|- |g \mathcal{W}_i \backslash \mathcal{W}_i|, \ \mathcal{W}_i:= \{ \text{hyperplanes separating $o$ from $Z_i$} \}$$
where $o \in X$ is a fixed basepoint. Finally, we prove that all the elements in the kernel of $\bigoplus_i \mathfrak{h}_i : H \to \mathbb{Z}^k$ are elliptic.
\end{itemize}

\paragraph{CAT(0) geometry may be inefficient.} Among the results coming from CAT(0) geometry that are sometimes applied to CAT(0) cube complexes, we can mention the structure of minimising sets of loxodromic isometries \cite[Proposition~II.6.2]{MR1744486} and the flat torus theorem \cite[Theorem~II.7.1]{MR1744486}. Most of the time, the cube complex is supposed to be finite-dimensional in order to avoid parabolic isometries (which may exist in infinite-dimensional CAT(0) cube complexes; see for instance \cite{ParaThompson, MR3198728} for a natural example (but parabolic isometries can be found in more elementary examples, such as the lamplighter group $\mathbb{Z} \wr \mathbb{Z}$)). However, this assumption is often not necessary in the results thus obtained. This is essentially due to the fact that, up to a subdivision, an isometry of a median graph either fixes a vertex (\emph{elliptic}) or acts as a translation on some bi-infinite geodesic (\emph{loxodromic}) \cite{Axis}. This dichotomy holds regardless of the dimension of the underlying cube complex. 

As an example, consider the mapping class group $\mathrm{MCG}_g$ of a closed surface of genus $g \geq 3$. It is known that $\mathrm{MCG}_g$ cannot act properly on a CAT(0) space by semisimple isometries \cite{MR1411351}, and, more precisely, that Dehn twists must be elliptic for every action of $\mathrm{MCG}_g$ on a CAT(0) space by semisimple isometries \cite{MR2665003}. Because there is no parabolic isometries in finite-dimensional CAT(0) cube complexes, one obtains severe restrictions on possible actions of $\mathrm{MCG}_g$ on finite-dimensional CAT(0) cube complexes\footnote{Actions of mapping class groups on CAT(0) cube complexes are of interest in particular due to connections to other well-known open questions: constructing actions on CAT(0) cube complexes with no global fixed points would imply that mapping class groups do not satisfy Kazhdan's property (T); and proving the fixed-point property on finite-dimensional CAT(0) cube complexes would imply that mapping class groups do not virtually surjects onto $\mathbb{Z}$.}. It turns out that, by following the same arguments but with a median perspective, the same conclusions can be obtained without any restriction on the dimension: Dehn twists are elliptic for every action of $\mathrm{MCG}_g$ on a median graph \cite{MR4574362}. The argument goes as follows.

Let $G$ be an arbitrary group acting on a median graph $X$. Assume that $G$ contains a central element $z$ having unbounded orbits in $X$. 
\begin{itemize}
	\item Up to subdividing $X$, we assume that $z$ admits an axis $\gamma$. For every $g\in G$, $g \gamma$ is an axis for $gzg^{-1}=z$. But different axes for the same element cross exactly the same hyperplanes, which implies that all the $g \gamma$ have a point at infinity in the same component $Y$ of the Roller boundary of $X$. This component must be stabilised by $G$.
	\item Then we get a morphism $\Theta : G \to \mathbb{Z}$ defined by $$g \mapsto |\mathcal{W}(o|Y) \backslash \mathcal{W}(go|Y)| - |\mathcal{W}(go|Y) \backslash \mathcal{W}(o|Y)|,$$ where $\mathcal{W}(\cdot | \cdot)$ denotes the set of the hyperplanes separating the two subsets under consideration and where $o \in X$ is a fixed basepoint. 
	\item By taking $o$ on $\gamma$, one sees that $\Theta(z)$ coincides with the translation length of $z$. In particular, $\Theta(z) \neq 0$.
\end{itemize}
Thus, we have proved that the central element $z$ survives in some quotient $G \twoheadrightarrow \mathbb{Z}$. Now, let $\mathrm{MCG}_g$ act on a median graph and let $\tau$ be a Dehn twist along a simple closed curve $c$. If $\tau$ has unbounded orbits in the median graph, then it follows from what we have just said that the centraliser $C(\tau)$ of $\tau$ must surject onto $\mathbb{Z}$ such that $\tau$ has a non-trivial image. But the mapping class group of the surface cut along $c$ surjects onto $C(\tau)$ and has finite abelianisation. So $\tau$ must have bounded orbits.

\medskip
Another central result of CAT(0) geometry is the fixed-point property. Namely, if a group acts on a complete CAT(0) space with bounded orbits, then it fixes some point globally. Because CAT(0) cube complexes with infinite cubes are not complete, the fixed-point property is sometimes stated for CAT(0) cube complexes of finite dimension only. However, regardless of the cubical dimension of a median graph $X$, every group $G$ that acts on $X$ with bounded orbits always stabilises a cube globally. In fact, this can be already deduced from the CAT(0) fixed-point property. Indeed, given a group $G$ acting on a CAT(0) cube complex $X$, we know that $G$ fixes a point in the completion of $X$. But this completion can be described as the closure $\bar{X}$ of the image of $X$ into $\ell^2(\mathcal{H})$ under the canonical embedding, where $\mathcal{H}$ denotes the set of the hyperplanes of $X$. And one can show that the barycentre of the bounded $G$-orbit of a vertex coming from $X$, which is fixed by $G$ but which belongs a priori only to $\bar{X}$, actually belongs to $X$. (See \cite{MR1459140} for the case of a cube.) However, the median approach is instructive \cite{Roller}. A first step reduces the problem to finite orbits, where the following argument can be applied:
\begin{itemize}
	\item If a group $G$ acts on a median graph $X$ with a finite orbit $\mathcal{O}$, then, up to replacing $X$ with the convex hull of $\mathcal{O}$, we can assume that $X$ is finite.
	\item Given a hyperplane $J$, say that $J$ is \emph{balanced} if its two halfspaces have the same cardinality. Otherwise, say that $J$ is \emph{unbalanced} and denote by $J^+$ its larger halfspace. Clearly, the $J^+$ pairwise intersect, so the Helly property implies that their global intersection $Q$ is non-empty.
	\item The hyperplanes of $Q$ correspond to the balanced hyperplanes of $X$. But the balanced hyperplanes of $X$ are clearly pairwise transverse, so $Q$ must be a cube. It is stabilised by $G$ by construction. 
\end{itemize}
This elementary argument has been applied successfully, in one form or another, in several situations. For instance, to amenable groups \cite{MR3833346} and to purely elliptic actions \cite{GLU}. 

\medskip
As another example of how median graphs can be used to obtain more efficient statements, consider the local criterion provided by Theorem~\ref{thm:Local}. Given a cube complex $X$ and a vertex $v \in X$, for all neighbours $v_1, \ldots, v_n \in X$ such that the edges $[v,v_1], \ldots, [v,v_n]$ pairwise span squares, we need to verify that $[v,v_1], \ldots, [v,v_n]$ span an $n$-cube. Therefore, we need to understand, at least a little bit, the ball centred at $v$ of radius $n$. In fact, from the point of view of median graphs, it suffices to investigate balls of radius $3$:

\begin{thm}\label{thm:VeryLocal}\footnote{The theorem is essentially contained in the proof of \cite[Theorem~6.1]{MR1748966}. The philosophy underlying the statement is probably well-known to specialists. For instance, the approach to CAT(0) cube complexes based on disc diagrams described in \cite{MR4298722} only uses the $2$-skeletons of cube complexes.}
Let $X$ be a (simplicial) graph. Then $X$ is median if and only if the following conditions hold:
\begin{itemize}
	\item the square-completion of $X$ (i.e.\ the square complex obtained from $X$ by filling in all the $4$-cycles) is simply connected;
	\item for every vertex $v \in X$ and all neighbours $a,b \in X$, the edges $[v,a],[v,b]$ span at most one $4$-cycle;
	\item for every vertex $v \in X$ and all neighbours $a,b,c \in X$, if the edges $[v,a],[v,b],[v,c]$ pairwise span a $4$-cycle then they globally span (the one-skeleton of) a $3$-cube.
\end{itemize}
\end{thm}

For explicit constructions of median graphs, this may simplify significantly the verification of the local criterion. (For instance, we applied this criterion in \cite{ARMCG2}.) This observation may also be useful in cubulating quotients, since the local criterion is preserved under quotients by subgroups with large injectivity radius, even if the cube complex is infinite-dimensional.

\paragraph{Median geometry is conceptually important.} Recently, median geometry and its related combinatorics of hyperplanes have been adapted in several directions. Variations around median geometry include coarse median spaces \cite{MR3037559}, as well as the subfamilies given by hierarchically hyperbolic spaces \cite{MR3650081, MR3956144} and locally quasi-cubical spaces \cite{LQC}. They explore similarities between median graphs and mapping class groups. The machinery of hyperplanes also found applications beyond the scope of median graphs, including CAT(0) spaces themselves \cite{HypModel, ZSurvey}. All these developments are fairly recent and require to be further investigated, but they highlight the conceptual potential of median geometry and its related concepts.

\paragraph{Historical background.} One could argue that it would be reasonable to keep the terminology ``CAT(0) cube complexes'' in order to keep track of the historical roots of the objects. But this would be a point of view restricted to geometric group theory, ignoring the perspective offered by metric graph theory. Median graphs have an even older history in metric graph theory. Highlighting this history would be helpful in order to promote interactions between geometric group theory and metric graph theory. Despite the fact that many concepts related to CAT(0) cube complexes in geometric group theory can be found in metric graph theory several years earlier\footnote{A notable example is provided by Buneman graphs of cut systems \cite{MR1002234}, which turn out to coincide with cubulations of wallspaces \cite{MR2059193, MR2197811}.}, very few geometric group theorists have some background in metric graph theory. 

As an illustration of the benefit that can follow from the perspective offered by metric graph theory, let us mention the \emph{Djokovi\'{c}-Winkler relation}, which sheds an interesting light on the notion of hyperplanes. Given a graph $X$, say that two edges $\{a,b\}$ and $\{x,y\}$ are in \emph{$\theta$-relation} if
$$d(a,x)+d(b,y) \neq d(a,y)+d(b,x).$$
Alternatively, the edge $\{a,b\}$ yields a partition of $X$ into three vertices: the vertices closer to $a$ than to $b$, the vertices closer to $b$ than to $a$, and the vertices at equal distance to $a$ and $b$. The edge $\{x,y\}$ is then in $\theta$-relation with $\{a,b\}$ whenever $x$ and $y$ belong to two distinct pieces of the partition. 

The Djokovi\'{c}-Winkler relation is reflexive and symmetric, but it may not be transitive. So we let $\theta^\ast$ denote its transitive closure. 

 The $\theta$-relation was introduced in \cite{MR314669}, and latter generalised in \cite{MR727925} for non-bipartite graphs, in order to characterise \emph{partial cubes}, i.e.\ isometrically embedded subgraphs in hypercubes. Then, \cite{MR776391} proved that the relation plays a fundamental role in representations of graphs as subgraphs of Cartesian products. More precisely, given a graph $X$, there is a canonical and optimal isometric embedding
$$X \hookrightarrow \prod\limits_{e \in E(X)/ \theta^\ast} X // e^c$$
where each $X // e^c$ denotes the graph obtained from $X$ by collapsing all the edges not in the class $e$. 

When applied to median graphs, $\theta^\ast$ can be described as the reflexive-transitive closure of the relation that identifies two edges whenever they are opposite in some $4$-cycle. In other words, we recover the standard notion of hyperplanes. And the embedding above coincides with the well-known embedding of a median graph inside a hypercube (whose will be, in addition, a retract). Hence a much deeper perspective on the notion of hyperplanes in CAT(0) cube complexes, which echoes some constructions of cubulations that implicitly use the Djokovi\'{c}-Winkler relation (e.g.\ Coxeter groups \cite{MR1983376}).

\paragraph{Median graphs or median cube complexes?} So far, our goal was to promote the use of ``median graphs'' instead of ``CAT(0) cube complexes''. However, other alternatives, not mentioned so far, exist. For instance, ``median cube complexes'', in order to keep the higher-dimensional cellular structure. Here, a \emph{median cube complex} would refer to a cube complex such that the length metric extending the $\ell^1$-metrics on each cube (whose sides, by convention, have length one) is median.

Even in the study of median graphs, cubes are fundamental. (It is understood that, when dealing with graphs, a cube always refers to a (sub)graph isomorphic to a product of pairs of adjacent edges (in other words, to the one-skeleton of a topological cube).) For instance:
\begin{itemize}
	\item In a median graph, every group of isometries with bounded orbits stabilises a cube.
	\item Given a median graph of finite cubical dimension, there is a natural bijection between the maximal cubes and the maximal collections of pairwise transverse hyperplanes.
	\item Median graphs coincide with retracts of infinite cubes.
\end{itemize}
Moreover, cube-completions of median graphs are necessary to study topological and cohomological properties of groups, such as finiteness properties (especially when combined with Morse theory \cite{MR1465330}). They are also useful to recognise or construct examples of groups satisfying fixed-point properties, because the Helly property can be applied to fixed-point sets in cube-completions (which are contractible). 

Another argument in favour of keeping actual cubes in our spaces comes from the fact that, as already mentioned, some group actions on median graphs are constructed by taking the universal covers of cube complexes satisfying the local criterion given by Theorem~\ref{thm:Local} (which could be called \emph{locally median cube complexes}). 

On the other hand, many other group actions are obtained by constructing the median graph directly. Then, adding the cubes would sound inefficient (for instance, the action would not be minimal, since the one-skeleton would provide an invariant median subspace). Moreover, as mentioned in \S 2, this is the distances between vertices that contain the interesting information. 

In conclusion, despite the fact that CAT(0) geometry is clearly not relevant in the study of CAT(0) cube complexes, it is not always desirable, and sometimes not even pertinent at all, to avoid the cube complex structure. This is why both median graphs and their cube-completions (which can be referred to as \emph{median cube complexes}) have to be considered.

\paragraph{CAT(0) geometry somewhere?} In order to justify the replacement of CAT(0) cube complexes with median graphs, we do not have to prove that CAT(0) geometry will never be relevant. Actually, there may (and probably do) exist a few specific cases where CAT(0) geometry on cube complexes is superior from some perspective. But it suffices to show, as tried in this note, that median geometry is more natural, efficient, and satisfying in most cases.

I would like to emphasize that the change of terminology suggested in this note does not aim at completely ignoring nonpositive curvature. Indeed, nonpositive curvature is a useful conceptual guide, so this is something to keep in mind. For instance, it helps guessing results to prove and it explains similarities with other geometries. But this is not something specific to CAT(0) geometry. In the literature, there are many geometries that are thought of as nonpositively curved but also quite different from CAT(0) geometry. 

Also, our comments here deal with geometric group theory only. In other contexts, the CAT(0) geometry may be the most relevant. For instance, this is the case in applications to robotics. Indeed, state spaces of some robots turn out to be often CAT(0) cube complexes (see for instance \cite{MR2301699}). Moving such a robot efficiently amounts to finding a geodesic between two states in the cube complex with respect to the CAT(0) metric. Considering the median-metric here would not be reasonable. For instance, given an articulated arm, moving one articulation and then another could define a median-geodesic, while a CAT(0)-geodesic would require to move the two articulations simultaneously (which corresponds to what we have in mind as an efficient move).

\addcontentsline{toc}{section}{References}

\bibliographystyle{alpha}
{\footnotesize\bibliography{CCvsMedian}}

\Address

%\addcontentsline{toc}{section}{Index}
%
%\printindex

\end{document}